%
%
%
%
%
%
%
%
%
%
%
%
%
%
\scrollmode
\magnification=\magstep1
\parskip=\smallskipamount

\hoffset=1cm \hsize=12cm

\def\demo#1:{\smallskip \it{#1}. \rm}
\def\ni{\noindent}               
\def\ll{\leftline}
\def\cl{\centerline}

%
%
\outer\def\beginsection#1\par{\bigskip
  \message{#1}\centerline{\bf #1}
  \nobreak\smallskip\vskip-\parskip\noindent}

%
%
\outer\def\proclaim#1:#2\par{\medbreak\vskip-\parskip
	{\rm#1.\enspace}{\sl#2}
  \ifdim\lastskip<\medskipamount \removelastskip\penalty55\smallskip\fi}

%
%

\def\R{{\rm I\kern-0.2em R\kern0.2em \kern-0.2em}}
\def\N{{\rm I\kern-0.2em N\kern0.2em \kern-0.2em}}
\def\P{{\rm I\kern-0.2em P\kern0.2em \kern-0.2em}}
\def\B{{\rm I\kern-0.2em B\kern0.2em \kern-0.2em}}
\def\C{{\rm C\kern-.4em {\vrule height1.4ex width.08em depth-.04ex}\;}}
\def\CP{\C\P}

%
%
%
%

\def\cC{{\cal C}}

\def\cF{{\cal F}}

\def\cO{{\cal O}}
\def\cP{{\cal P}}

\def\cS{{\cal S}}

%
%
%
\def\a{\alpha}
\def\b{\beta}
\def\g{\gamma}
\def\d{\delta}
\def\e{\epsilon}
\def\z{\zeta}


%
%
%
%
\def\bar{\overline}              
\def\bs{\backslash}              

\def\di{\partial}                
\def\dibar{\bar\partial}         

%
%
\def\dim{{\rm dim}\,}                    
\def\holo{holomorphic}                   
\def\nbd{neighborhood}                   
\def\spsc{strongly\ pseudoconvex}        
\def\psh{plurisubharmonic}               
\def\spsh{strongly\ plurisubharmonic}
\def\tr{totally real}                    
\def\hc{holomorphically convex}          
\def\ss{\subset\!\subset}                
\def\supp{{\rm supp}\,}                  
\def\iff{if and only if}

\def\hra{\hookrightarrow}


\def\wt{\widetilde}


\def\begin{\ll{}\vskip 10mm \nopagenumbers}  
\def\pn{\footline={\hss\tenrm\folio\hss}}   

%
%
%
%

\begin	
\cl{\bf HOLOMORPHIC SUBMERSIONS FROM}
\smallskip
\cl{\bf STEIN MANIFOLDS}
\vskip 8mm
\centerline{\bf by Franc FORSTNERI\v C}
\bigskip

\vskip 7mm
\cl{\bf --------------}
\vskip 7mm

\cl{\bf 1. Introduction.} 
\medskip
\rm

Recall that a {\it Stein manifold} is a complex manifold biholomorphic 
to a closed complex submanifold of a complex Euclidean space $\C^N$
[GR, p.\ 226]. A holomorphic map $f\colon X\to Y$ whose
differential $df_x\colon T_x X\to T_{f(x)} Y$ is surjective 
for every $x\in X$ is said to be a {\it submersion} of $X$ to $Y$. 
The following was proved in [F1, Theorem II]; for $n=q=1$ see also [GN]:

{\it A Stein manifold $X$ admits a \holo\ submersion to $\C^q$ 
for some $q< \dim X$ \iff\ its tangent bundle $TX$ admits a surjective 
complex vector bundle map onto the trivial rank $q$ bundle $\tau_X^q=X\times \C^q$.}

The necessity of the latter condition is clear since the tangent map 
of a submersion $X\to\C^q$ induces a surjective vector bundle map 
$TX\to\tau_X^q$. The corresponding problem for $n=q>1$ remains open [F1]. 

In this paper we consider the analogous problem 
with the target space $\C^q$ replaced by a more general complex 
manifold $Y$. The following result in the smooth category
was proved by A.\ Phillips [P] and M.\ Gromov [Gr1, Gr3]:

\smallskip
{\it A continuous map $f_0\colon X\to Y$ from a smooth open manifold $X$ 
to a smooth manifold $Y$ is homotopic to a 
submersion $f_1\colon X \to Y$ if and only if there exists 
a surjective vector bundle map $\iota_0\colon TX\to f_0^*TY$ over $X$,
i.e., the pull-back $f_0^* TY$ is a quotient bundle of $TX$.  
The regular homotopy classes of submersions $X\to Y$ 
are classified by the homotopy classes of their tangent maps.} 
\smallskip

In the holomorphic category the topological condition on the 
existence of $\iota_0$  does not suffice.  
For example, a Kobayashi-hyperbolic manifold $Y$ admits no nonconstant 
holomorphic images of $\C^n$ but the topological condition trivially 
holds for the constant map $\C^n\to y_0 \in Y$. This suggests that 
we restrict the attention to manifolds which admit 
`sufficiently many' holomorphic submersions $\C^n\to Y$ for a given 
$n \ge \dim Y$. It turns out that a suitably precise form of this 
condition suffices for a complete analogue of the Gromov-Phillips 
theorem, thus justifying the following heuristic principle
(a form of the Oka principle):

%
%
\smallskip
{\it If a complex manifold $Y$ admits sufficiently many holomorphic 
submersions $\C^n\to Y$ for a given integer $n\ge \dim Y$ 
then the existence of holomorphic submersions $X\to Y$ from $n$-dimensional 
Stein manifolds $X$ reduces to a purely topological problem.}
\smallskip

We say that a complex manifold $Y$ satisfies 
{\it Property $S_n$} for some integer $n\ge \dim Y$ if any holomorphic 
submersion $O\to Y$ from a compact convex set $O\subset \C^n$ of a 
certain special type can be approximated uniformly on $O$ by 
holomorphic submersions $\C^n\to Y$ (Definition 1 in Sect.\ 2). 
The following is a special case of our main result, Theorem 2.1:

\pn

\proclaim THEOREM: 
Assume that $X$ is an $n$-dimensional Stein manifold and $Y$ 
is a complex manifold satisfying Property $S_n$. 
A continuous map $f\colon X\to Y$ is homotopic to a 
\holo\ submersion of $X$ to $Y$ \iff\ there exists a surjective 
complex vector bundle map $\iota\colon TX\to f^* TY$.

Furthermore, {\it Property $S_n$ is both necessary and sufficient
for a stronger version of the above theorem with approximation
on compact $\cO(X)$-convex subsets of the source 
manifold $X$}; Theorem 2.1 (b).
If $Y$ also satisfies the analogous {\it Property $HS_n$} 
concerning the approximation of homotopies of holomorphic
submersions $O\to Y$ by homotopies of submersions $\C^n\to Y$ 
then for any $n$-dimensional Stein manifold $X$ the regular 
homotopy classes of holomorphic submersions $X\to Y$ are 
classified by the homotopy classes of their tangent maps 
(Corollary 2.4).

If $n=\dim X \ge 2\, \dim Y-1$ then $f^* TY$ is always a 
quotient of $TX$ by topological reasons and hence every continuous 
map $X\to Y$ is homotopic to a \holo\ submersion (Corollary 2.3).

One cannot use Property $S_n$ directly since a general submersion 
$X\to Y$ does not factor as $X\to\C^n\to Y$. Instead such decompositions 
are used on small subsets of $X$ and the resulting local submersions 
of $X$ to $Y$ are pieced together into a global submersion by 
the analytic tools developed in [F1] and in this paper.

Properties $S_n$ and $HS_n$ are equivalent to apparently 
weaker conditions on uniform approximability of submersions 
$O\to Y$ on special compact convex sets $O\subset \C^n$
by submersions $\C^n\bs A\to Y$, where $A\subset\C^n$ is an algebraic 
subvariety of codimension at least two which does not intersect $O$
(this relies upon the theory of holomorphic automorphisms of $\C^n$). 
%
%
In Sect.\ 5 we establish Properties $S_n$ and $HS_n$ when $n>\dim Y=q$ 
(or $n=q=1$) and $Y$ is any of the following:
$\C^q$, $\CP^q$, a complex Grassmanian, a Zariski open set with thin 
complement (containing no hypersurfaces) in any of the above, or
a holomorphic quotient of any of the above (this class contains 
all complex tori and Hopf manifolds).

For Riemann surfaces we obtain a complete answer
by proving that the following are equivalent (Corollary 2.8):
\item{(a)} $Y$ is one of the Riemann surfaces $\CP^1$, $\C$, 
$\C^*=\C\bs \{0\}$, or a complex torus (the quotient $\C/\Gamma$
by a rank two lattice $\Gamma\subset \C$);
\item{(b)} any continuous map from a Stein manifold $X$ to $Y$ 
is homotopic to a holomorphic submersion.

The Riemann surfaces listed in (a) are precisely those which are not Kobayashi 
hyperbolic. When $Y$ is $\C^*$ or a torus our result is new even 
when $X$ is an open Riemann surface.

The heuristic principle behind our main result is reminiscent of 
Gromov's extension of the Oka-Grauert principle [Gr4]: 
{\it The existence of many dominating \holo\ maps $\C^n\to Y$
implies the existence of many holomorphic maps $X\to Y$ from any Stein 
manifold $X$}. A dominating map $s\colon \C^n\to Y$ based at a point $y \in  Y$
is one for which $s(0)=y$ and $ds_0\colon T_0\C^n \to T_y Y$
is surjective (i.e., $s$ is a submersion near $0\in \C^n$).
A {\it dominating spray} on $Y$ is a family of such maps depending holomorphically 
on $y\in Y$ (its  domain is a holomorphic vector bundle over $Y$). 
Comparing with Property $S_n$ we see that the two conditions have 
similar flavour. For the Oka-Grauert-Gromov theory we refer to 
[G3, Gr4, HL2, FP1, FP2, FP3, F3, W].

An important ingredient in our construction of submersions 
is a holomorphic approximation theorem on certain handlebodies 
in arbitrary complex manifolds; Theorem 3.2 in Sect.\ 3.
This result, together with a geometric lemma from [F1],
gives an approximate extension of the submersion across a 
critical level of a \spsh\ exhaustion function
in the source manifold. In Sect.\ 6 we use the same method 
to give a simple proof of the Oka-Grauert principle for sections 
of holomorphic fiber bundles over Stein manifolds whose fiber 
admits a finite dominating family of sprays.

\beginsection 2. The main results.

We denote by $\cO(X)$ the algebra of all holomorphic functions on
a complex manifold $X$.  A function (or map) is 
holomorphic on a closed subset $K$ in $X$ if it is holomorphic in 
some open \nbd\ of $K$; the set of all such functions (with the usual 
identification of functions which agree in a \nbd\ of $K$) will be denoted $\cO(K)$. 
Any statement concerning a holomorphic map on a closed set should 
be understood in the sense that it holds in some open \nbd;
for homotopies of maps the \nbd\ is the same for all maps in 
the homotopy. A compact set $K \subset X$ is $\cO(X)$-convex
if for every $p\in X\bs K$ there exists $f\in \cO(X)$ such that
$|f(p)| >\sup_K |f|$.

A {\it homotopy of holomorphic submersions} $X\to Y$ is a family of 
holomorphic submersions $f_t\colon X\to Y$ $(t\in [0,1])$ depending 
continuously on $t$. It follows that the tangent maps $Tf_t \colon TX\to TY$
are also continuous in $t$.

Let $z=(z_1,\ldots,z_n)$ be the coordinates on $\C^n$, with 
$z_j=x_j+iy_j$. Set
$$
	Q =\{z\in \C^n\colon |x_j| \le 1,\ |y_j|\le 1,\ j=1,\ldots,n\}.  \eqno(2.1)
$$
A {\it special convex set} in $\C^n$ is a compact convex subset of the form
$$ 
	O=\{z\in Q \colon y_n \le h(z_1,\ldots,z_{n-1},x_n)\},   \eqno(2.2)
$$
where $h$ is a smooth (weakly) concave function with values in $(-1,1)$.

%
%
\proclaim DEFINITION 1: 
Let $d$ be a distance function induced by a smooth Riemannian metric
on complex manifold $Y$.
\hfill\break
(a) $Y$ satisfies {\rm Property} $S_n$ if for any holomorphic submersion 
$f\colon O\to Y$ on a special convex set $O\subset \C^n$ 
and any $\e>0$ there is a  \holo\ submersions $\wt f\colon Q\to Y$
satisfying $\sup_{x\in O} d(f(x),\wt f(x)) <\e$.
\hfill \break  
(b) $Y$ satisfies {\rm Property} $HS_n$ if for any homotopy of 
holomorphic submersions $f_t\colon O\to Y$ $(t\in [0,1])$ 
such that $f_0$ and $f_1$ extend to holomorphic submersions 
$Q\to Y$ there exists for any $\e>0$ a homotopy of holomorphic 
submersions $\wt f_t\colon Q \to Y$ $(t\in [0,1])$ satisfying
$\wt f_0=f_0$, $\wt f_1=f_1$, and  
$\sup_{x\in O,\ t\in [0,1]} d(f_t(x),\wt f_t(x)) <\e$.
\medskip

An obvious induction shows that Property $S_n$ of $Y$
implies the following: 
{\it Any  holomorphic submersion 
$f\colon O\to Y$ on a special convex set $O\subset \C^n$
can be approximated uniformly on $O$ by holomorphic submersions
$\C^n\to Y$.} The analogous remark holds for Property $HS_n$.

Let $X$ and $Y$ be complex manifolds. Denote by $\cS(X,Y)$ 
the set of all pairs $(f,\iota)$ where $f\colon X\to Y$ 
is a continuous map and $\iota\colon TX\to TY$ is a fiberwise 
surjective complex vector bundle map making the following diagram 
commute:  
$$ \matrix{    & {}\atop \iota \cr 
            TX & \longrightarrow & TY \cr
             \downarrow & {}\atop f & \downarrow \cr
             X & \longrightarrow & Y  \cr}
$$
\cl{\it Figure 1: The set $\cS(X,Y)$} 
\vskip 5mm

Such $\iota$ is the composition of a surjective complex vector bundle 
map $TX\to f^* TY$ with the natural map $f^* TY\to TY$.
Let $\cS_{holo}(X,Y)$ consist of all pairs $(f,Tf)\in  \cS(X,Y)$ 
where $f\colon X\to Y$ is a holomorphic submersion and $Tf$ its 
tangent map. We equip $\cS(X,Y)$ with the compact-open topology. 

The following is our main theorem.

\proclaim THEOREM 2.1:
Assume that $X$ is a Stein manifold of dimension $n$ and $Y$ 
is a complex manifold satisfying Property $S_n$. 
\item{(a)} {\rm (Existence of submersions)} 
Every $(f_0,\iota_0) \in\cS(X,Y)$ can be connected by a path 
$\{(f_t,\iota_t)\}_{t\in [0,1]} \subset \cS(X,Y)$ to some 
$(f_1,Tf_1) \in  \cS_{holo}(X,Y)$. 
\item{(b)} {\rm (Approximation)} 
If $K\subset X$ is a compact $\cO(X)$-convex subset, 
$f_0|_K$ is a holomorphic submersion, and $\iota_0|_K = Tf_0|_K$ 
then the path $(f_t,\iota_t)$ in (a) can be chosen 
such that for every  $t\in [0,1]$, $f_t|_K$ is a holomorphic
submersion uniformly close to $f_0|_K$ and $\iota_t|_K=T f_t|_K$.
\item{(c)} {\rm (Regular homotopies of submersions)}
If $Y$ satisfies Property $HS_n$ then any path 
$\xi_t=(f_t,\iota_t) \in \cS(X,Y)$ 
$(t\in [0,1])$ with $\xi_0,\ \xi_1\in \cS_{holo}(X,Y)$ can be  
deformed with fixed ends to a path in $\cS_{holo}(X,Y)$.
\medskip

Theorem 2.1 is proved in Section 4. In the special case $Y=\C^q$, 
with $n>q$ or $n=q=1$, parts (a) and (b) were proved in [F1]
(for $n=q=1$ see also [GN]), but part (c) is new even in this case. 
The relevant Property $HS_n$ holds for $\C^q$ whenever $n>q$ 
(Proposition 2.5 below).

\smallskip
\it Remarks. \rm 
1.\ If the conclusion of part (b) in Theorem 2.1 holds for a given complex 
manifold $Y$, with $K$ any compact convex set in $X=\C^n$, then (by definition) 
$Y$ satisfies Property $S_n$. Hence both the topological condition 
(the existence of $\iota$) and the analytic condition 
(Property $S_n$ of $Y$) are necessary for parts (a) and (b) 
in Theorem 2.1. 

2.\ Clearly $HS_n \Rightarrow S_n$. We don't know whether the 
converse always holds, but in all examples for which we prove $S_n$ 
we also prove $HS_n$. 
\smallskip

Part (a) of Theorem 2.1 implies the following.

\proclaim COROLLARY 2.2: 
If $X$ is an $n$-dimensional Stein manifold and  $Y$ satisfies 
Property $S_n$ then for any $(f,\iota)\in \cS(X,Y)$ there is 
a nonsingular holomorphic foliation of $X$ whose normal bundle 
is isomorphic to $f^*TY$ (as a topological complex vector bundle
over $X$).

Such a foliation is given by the level sets of a submersion $X\to Y$ 
furnished by Theorem 2.1 (a). The corresponding result for foliations
with trivial normal bundle (and $Y=\C^q$) was obtained in [F1].

If $\dim X \ge 2\dim Y-1$ then every map $f\colon X\to Y$ can be covered
by a map $\iota \colon TX\to TY$ such that $(f,\iota)\in\cS(X,Y)$. 
This follows by standard topological methods from the fact 
that an $n$-dimensional Stein manifold is homotopic to an $n$-dimensional
CW-complex; it will also be clear from our proof of 
Theorem 2.1. Hence Theorem 2.1 implies

\proclaim COROLLARY 2.3: If $Y$ satisfies Property $S_n$ 
for some $n\ge 2\, \dim Y - 1$ then any map $f\colon X\to Y$ from an 
$n$-dimensional Stein manifold $X$ is homotopic to a 
holomorphic submersion of $X$ to $Y$.

Combining parts (a) and (c) in Theorem 2.1 we also obtain

\proclaim COROLLARY 2.4: 
If a complex manifold $Y$ satisfies Property $HS_n$
then for every $n$-dimensional Stein manifold $X$ the natural inclusion 
$$
	\cS_{holo}(X,Y)\hra \cS(X,Y)
$$ 
induces a bijective correspondence of the path-connected components of the two spaces.
\smallskip

Comparing with the Oka-Grauert principle [G3, Gr4, FP1] one might 
expect that the above inclusion is a {\it weak homotopy equivalence}; 
however, our proof does not give this much.

In the remainder of this section we discuss the question
which manifolds satisfy Properties $S_n$ and $HS_n$.
For Riemann surfaces both properties hold precisely on 
the non-hyperbolic ones (Corollary 2.8).
For manifolds of dimension $>1$ a complete answer 
seems out of reach even in the class of projective 
algebraic surfaces. The following result gives some 
nontrivial examples.

\proclaim PROPOSITION 2.5: The following manifolds $Y^q$ satisfy 
Properties $S_n$ and $HS_n$ for any $n>q$, as well as for $n=q=1$:
\item{(a)} $\C^q$, $\CP^q$, or a complex Grassmanian;
\item{(b)} a Zariski open subset in any of the manifolds 
from (a) whose complement contains no complex hypersurfaces.
\smallskip

For a proof see Sect.\ 5. The following simple observation will be useful.

\proclaim PROPOSITION 2.6: Let $\pi \colon \wt Y\to Y$ be a \holo\ covering. 
If one of the manifolds $Y$, $\wt Y$ satisfies Property 
$S_n$ (resp.\ $HS_n$) then so does the other.

\demo Proof: Assume first that $\wt Y$ satisfies $S_n$.
Given a \holo\ submersion $f\colon O\to Y$ from a special 
compact convex set $O\subset \C^n$,
there is a holomorphic lifting $g\colon O\to \wt Y$ (satisfying
$\pi\circ g=f$) which is also a submersion since $\pi$ is 
locally biholomorphic. By the assumed Property $S_n$ of $\wt Y$
we can approximate $g$ by an entire submersion $\wt g\colon \C^n\to \wt Y$,
and then $\wt f = \pi\circ \wt g \colon \C^n\to Y$ is a submersion
such that $\wt f|_O$ approximates $f$.  This proves that 
$Y$ also satisfies $S_n$. 

Conversely, assume that $Y$ satisfies $S_n$. Given a submersion
$g\colon O\to\wt Y$, we approximate the submersion $f=\pi\circ g\colon O\to Y$ 
by a submersion $\wt f \colon \C^n \to Y$ (using Property $S_n$ of $Y$)
and then lift $\wt f$ to a (unique!) map $\wt g\colon \C^n\to \wt Y$ 
satisfying $\wt g(z_0)= g(z_0)$ for some point $z_0\in O$. Then 
$\wt g$ is a \holo\ submersion which approximates $g$ on $O$. 

The same arguments hold for homotopies of submersions which 
gives the corresponding statement for Property $HS_n$. 
This proves Proposition 2.6.
\smallskip

Combining Propositions 2.5 and 2.6 we obtain

\proclaim COROLLARY 2.7:  A complex manifold $Y$ whose universal covering 
space is biholomorphic to $\C^q$, or to $\C^q\bs A$ for an algebraic
subvariety of codimension at least two, satisfies Properties
$S_n$ and $HS_n$ when $n>q$ or $n\ge q=1$.

The manifolds covered by $\C^q$ are of the form  
$\C^k\times (\C^*)^l \times T^s$ ($k+l+s=q$), where $\C^*=\C\bs \{0\}$ 
and $T^s$ is a complex torus, i.e., $T^s=\C^s/\Gamma$ for a lattice
$\Gamma\subset \C^s$ of rank $2s$. Recall that every {\it Hopf manifold}
is a holomorphic quotient of $\C^q\bs \{0\}$ and hence Theorem
2.1 applies (see [B, p.\ 172]).

The following result, which follows from
Proposition 2.5 and Corollary 2.7, is a complete solution to 
the submersion problem if $Y$ is a Riemann surface.

\proclaim COROLLARY 2.8: If $Y$ is any of the Riemann surfaces 
$\CP^1$, $\C$, $\C^*$, or a complex torus $\C/\Gamma$ then any 
continuous map $f_0\colon X\to Y$ from a Stein manifold $X$ is 
homotopic to a holomorphic submersion $f_1\colon X \to Y$. Furthermore,
if $f_0|_K\colon K\to Y$ is a holomorphic submersion for a compact
$\cO(X)$-convex subset $K\subset X$ then the homotopy can be chosen to 
appoximate $f_0$ uniformly on $K$. 
Conversely, if a Riemann surface $Y$ admits a nonconstant 
\holo\ map $\C\to Y$ then $Y$ belongs to the above list.

Indeed, the universal covering of any Riemann surface is either 
$\CP^1$, $\C$ or the disc, and only the first two admit a nonconstant 
holomorphic image of $\C$; their holomorphic quotients are 
listed in Corollary 2.8.

The manifolds for which we prove Property $S_n$ in this paper are all  
{\it subelliptic} in the sense of [F3], i.e., they admit a finite 
dominating family of holomorphic sprays; hence by the main result in
[F3] any continuous map $X\to Y$ is homotopic to a holomorphic map 
(the Oka-Grauert principle).

\smallskip
\it Problem 1. \rm
Does every subelliptic manifold $Y$ satisfy Property $S_n$ for 
all $n>\dim Y$~?
\smallskip

A good test case may be complex Lie groups $G$; 
a dominating holomorphic spray on $G$ is given by $s(g,v)=exp(v) g$ 
where $g\in G$ and $v$ is a vector in the Lie algebra of $G$.

\smallskip \it Problem 2. \rm 
Does a manifold $Y$ which satisfies $S_n$ for some $n\ge \dim Y$
also satisfy $S_k$ for $k>n$~?

\smallskip
\it Problem 3. \rm
Is the Property $S_n$ invariant with respect to proper 
modifications such as blow-ups and blow-downs~?

%
%
%
%
\beginsection 3. Holomorphic approximation on handlebodies.

The main result of this section, Theorem 3.2, concerns \holo\ approximation of 
mappings between complex manifolds or, more generally, of sections
of holomorphic submersions, on certain compact sets obtained by attaching 
a \tr\ submanifold to a compact holomorphically convex subset. 
Theorem 3.2, which generalizes a result of
of H\"ormander and Wermer [HW, Theorem 4.1], will be used in the proof 
of Theorem 2.1 in Sect.\ 4. The proof uses H\"orman\-der's 
$L^2$-solution to the $\dibar$-equation [H\"o1, H\"o2] on special 
Stein neighborhoods of $S$ furnished by Theorem 3.1 below. 
For $X=\C^n$ Theorem 3.1 is due to H\"ormander and Wermer [HW, Theorem 3.1]; 
here we extend it to arbitrary complex manifolds. 

A compact set $K$ in a complex manifold $X$ is said to be 
{\it holo\-mor\-phically convex} if $K$ has an open  Stein \nbd\ 
$\Omega$ in $X$ such that $K$ is $\cO(\Omega)$-convex.
By the classical theory (see e.g.\ Chapter 2 in [H\"o2]) 
holomorphic convexity of $K$ is equivalent to the existence of a Stein \nbd\ 
$\Omega$ of $K$ and a continuous (or smooth) plurisubharmonic function 
$\rho\ge 0$ on $\Omega$ such that $\rho^{-1}(0)=K$ and  $\rho$ is strongly 
plurisubharmonic on $\Omega\bs K$. We may take $\Omega = \{\rho<c_1\}$ for 
some $c_1>0$; for any $c\in (0, c_1)$ the sublevel set $\{\rho<c\} \ss \Omega$ is 
then Stein and Runge in $\Omega$ (Sect.\ 4.3 in [H\"o2]).

A $\cC^1$ submanifold $M$ of a complex manifold $X$
is {\it \tr\/} in $X$ if for each $p\in M$ the tangent space
$T_p M \subset T_p X$ contains no complex line of $T_p X$. 

%
%
%
%
\proclaim THEOREM 3.1: Let $X$ be a complex manifold.
Assume that $S=K_0\cup M$ is a closed subset of an open set in $X$
such that $K_0$ is compact holomorphically convex and $M=S\bs K_0$ is a 
$\cC^1$ totally real submanifold. If there exists a compact 
holomorphically convex set $K_1 \subset S$ which is a relative \nbd\ of $K_0$ 
in $S$ then every compact set $K$ with $K_0\subset K\subset S$ is 
holomorphically convex. Fix any such $K$. Let $d$ be a distance 
function on $X$ induced by a smooth Riemannian metric. Fix a 
\nbd\ $N$ of $K_0$ in $X$. There are a constant $C>1$ and for any 
sufficiently small $\e>0$ a Stein domain $\omega_\e \subset X$ such that
\item{(i)} $\omega_\e$ contains all points of distance $<\e$ from $K$, 
\item{(ii)} all points $x\in\omega_\e\bs N$ have distance $< C\e$ from $M$.
\medskip

For $X=\C^n$ this is Theorem 3.1 in [HW] where the result is proved with $C=2$. 
In our application (Sect.\ 4) $K_0$ is the closure of a strongly pseudoconvex 
domain and $M$ is attached to $K_0$ along a submanifold 
$bM \subset bK_0=:\Sigma$ which is complex tangential 
in $\Sigma$, i.e., $T_p bM \subset  T^\C_p \Sigma = T_p \Sigma\cap JT_p\Sigma$
for every $p\in bM$, where $J$ is the almost complex structure operator on $TX$.
For such $K_0$ and $M$, and with $X=\C^n$, results on holomorphic convexity 
of $K_0\cup M$ can be found in [E] (Lemmas 3.3.1.\ and 3.4.3.), 
[Ro, Lemma 2], [FK], and other papers.

\demo Proof of Theorem 3.1:
By the assumption on $K_0$ there is a smooth \psh\ function 
$\rho_0\ge 0$ in an open Stein \nbd\ $U_0\subset X$ of $K_0$ such that 
$\rho_0^{-1}(0)=K_0$. We may assume that the neighborhood $N$ 
of $K_0$ in Theorem 3.1 is chosen such that $N\subset U_0$ and $N\cap S\subset K_1$.
Choose a sufficiently small $c_2>0$ such that $\{\rho_0 < c_2\}\ss N$. 
Also choose constants $0<c_0<c'_0<c_1<c'_1<c_2$.

Since Theorem 3.1 only concerns compact subsets of $S$, 
we may assume that $S$ is compact. Since $M$ is totally real, 
there is a $\cC^2$ strongly plurisubharmonic
function $\tau\ge 0$ in an open set $V_0 \supset M$ which vanishes 
quadratically on $M$ [HW]. Replacing $\tau$ by $c\tau$ for a suitable
$c>0$ and shrinking $V_0$ around $M$ if necessary we may assume that 
$$
	\tau(x)< d(x,M)^2, \qquad x\in V_0.                           \eqno(3.1)
$$
Since $\tau$ vanishes quadratically on $M$, there is a $C>0$ such that 
$$
	\tau(x)> d(x,M)^2/C^2, 
			\qquad x\in V_0\cap \{\rho_0\ge c_0\}.        \eqno(3.2)
$$

Recall that $K_1$ is assumed to be holomorphically convex.
Choose a smooth \psh\ function $\rho_1\ge 0$ in an open \nbd\ $U_1$ of 
$K_1$ such that $\rho_1^{-1}(0)=K_1$ and $\rho_1$ is \spsh\ 
off $K_1$. Since $\rho_1$ vanishes at least to second order at 
each point $x\in K_1$, we may assume (after rescaling 
$\rho_1$ if necessary) that 
$$ 
	\rho_1(x) < d(x,S)^2/C^2, \qquad x \in U_1 \cap \{\rho\le c_2\} 
							              \eqno(3.3)
$$
where $C$ is the constant from (3.2).
Choose a smooth cut-off function $\chi\ge 0$ on $X$ which equals one
on $\{\rho_0\le c_0\}$ and satisfies $\supp \chi \ss \{\rho_0 < c'_0\}$.
For a sufficiently small $\d>0$ the function $\tau_\d=\tau-\delta\chi$
is \spsh\ on $V_0 \cap\{\rho_0 \ge c_0\}$ (we may need to shrink $V_0$
around $M$). Clearly $\tau_\d=\tau$ on $V_0\cap \{\rho_0 \ge c'_0\}$.
Fix such a $\d$. Choose an open \nbd\ $V \subset U_1\cup V_0$ of $S$
and let $\rho\colon V\to\R_+$ be defined by 
$$
	\rho= \cases{ \rho_1       & on $V\cap \{\rho_0 < c_0\}$, \cr
	    \max\{\rho_1,\tau_\d\} & on $V\cap \{c_0 \le \rho_0\le c'_0\}$, \cr
	              \tau         & on $V\cap \{\rho_0 > c'_0\}$.} 
$$
It is easily verified that these choices are compatible provided that 
the \nbd\ $V$ of $S$ is chosen sufficiently small. (When checking 
the compatibility near $M\cap \{\rho_0=c'_0\}$ the reader should observe that, 
by (3.2) and (3.3), we have $\rho_1 < \tau$ and hence $\rho=\tau$ there. 
Near $M\cap \{\rho_0=c_0\}$ we have $\rho_1\ge 0$ while $\tau_\d<0$, hence
$\rho=\rho_1$.) The function $\rho\ge 0$ is plurisubharmonic, \spsh\ on 
$V\cap \{\rho_0\ge c'_0\}$ (where it equals $\tau$), and $\rho^{-1}(0)=S$.
For every sufficiently small $\e>0$ the set
$$ 
	\omega_\e = \{x\in V\colon \rho(x)<\e^2\} \ss V
$$
is a pseudoconvex open \nbd\ of $S$ which satisfies
$$ 
    \{x\in X\colon d(x,S)<\e\} \subset \omega_\e,\qquad
    \omega_\e \bs N \subset \{x \in X\colon d(x,M)< C\e\}.
$$
The first inclusion follows from $\rho(x)<d(x,S)^2$
which is a consequence of (3.1) and (3.3). The second inclusion 
is a consequence of (3.2) and the fact that $\rho=\tau$ on 
$V\bs N$.

It remains to show that the sets $\omega_\e$ are Stein.
Fix $\e$ and choose an increasing convex function 
$h_\e\colon (-\infty,\e^2)\to \R$ with $\lim_{t\to \e^2} h(t)=+\infty$. 
Then $h_\e\circ\rho$ is a \psh\ exhaustion function on $\omega_\e$. 
In order to obtain a {\it \spsh\ exhaustion function} on $\omega_\e$ we 
proceed as follows. Choose a smooth \spsh\ function $\xi \colon U_0\to \R$ 
(such $\xi$ exists since $U_0$ is Stein). Also choose a smooth 
cut-off function $\chi\ge 0$ on $X$ such that $\chi=1$ 
on $\{\rho_0\le c_1\}$ and 
$\supp\chi \ss \{\rho_0 < c'_1\}$. If $\d>0$ is chosen sufficiently small 
then $\wt\rho= \rho+\d \chi \xi$ is \spsh\ on $V$. 
Indeed, on $V\cap \{\rho_0 \le c_1\}$ we have
$\wt \rho=\rho+ \d \xi$ which is \spsh\ for every $\d>0$; 
on $V\cap \{\rho_0 > c_1\}$ the function $\rho$ is \spsh\ and 
hence so is $\wt \rho$ provided that $\d$ is chosen sufficiently small. 
For such $\d$ the function $h_\e\circ\rho + \wt\rho$ is a \spsh\ exhaustion 
on $\omega_\e$ and hence $\omega_\e$ is Stein.

This gives a desired Stein \nbd\ basis $\omega_\e$ of $S$ satisfying Theorem 3.1. 
The same proof applies to any compact subset $K\subset S$ containing $K_0$.
Alternatively one can apply the above proof with $\rho$ 
replaced by $\rho+\tau_K$ where $\tau_K\ge 0$ is a smooth function 
which vanishes to order $>2$ on $K$. This completes the proof
of Theorem 3.1.

\smallskip \it Remark. \rm
Since the function $\rho$ constructed above is \psh\ 
on the Stein manifold $\omega_{\e_0}=\{\rho<\e_0\}$ 
for some small $\e_0>0$, its sublevel sets 
$\omega_\e=\{\rho<\e\}$ for $\e \in (0,\e_0)$ 
are Runge in $\omega_{\e_0}$.

\proclaim THEOREM 3.2: 
Let $K_0$ and $S=K_0\cup M$ be compact \hc\ subsets in a complex
manifold $X$ such that $M=S\bs K_0$ is a \tr\ $m$-dimensional submanifold 
of class $\cC^r$. Assume that $r\ge m/2+1$ and 
let $k$ be an integer satisfying $0\le k\le r - m/2 -1$. 
Given an open set $U\subset X$ containing $K_0$ and 
a map $f\colon U\cup M\to Y$ to a complex manifold $Y$ 
such that $f|_U$ is holomorphic and $f|_M \in \cC^r(M)$, 
there exist open sets $V_j\subset  X$
containing $S$ and holomorphic maps $f_j\colon V_j\to Y$ 
$(j=1,2,3,\ldots)$ such that, as $j\to\infty$, the sequence
$f_j$ converges to $f$ uniformly on $K_0$ and in the $\cC^k$-sense 
on $M$. If in addition $X_0$ is a closed complex subvariety of $X$
which does not intersect $M$ and $s\in \N$ then we can choose the 
approximating sequence such that $f_j$ agrees to order $s$ with $f$ 
along $X_0\cap V$ for all $j=1,2,3,\ldots$. The analogous result holds for 
sections $f\colon X\to Z$ of any \holo\ submersion $h\colon Z\to X$.

The domains $V_j$ of $f_j$ in Theorem 3.2 may shrink to $S$
as $j\to\infty$.
Although none of the manifolds $X,Y,Z$ in Theorem 3.2 is assumed to be Stein,
we shall reduce the proof to that case. In our applications
$K_0$ will be a sublevel set $\{\rho\le c\}$ of a smooth
\spsh\ function on $X$ and $M$ will be a smooth \tr\ handle attached to the 
hypersurface $bK_0=\{\rho=c\}$ along a legendrian (complex tangential) 
sphere.

\demo Proof:
When $X=\C^n$, $Y=\C$ (i.e., $f$ is a function) and $k=0$, 
Theorem 3.2 coincides with Theorem 4.1 in [HW, p.\ 11]. 
For later purposes we recall the sketch of proof which uses 
the Stein \nbd\ basis $\omega_\e \subset X$ of $S$ 
furnished by Theorem 3.1 in [HW] (compare with Theorem 3.1 above). 
One first obtains a $\cC^r$ extension $u$ of $f$ in a \nbd\ of $M$ 
such that $|D^s(\dibar u)|=o(d_M^{r-1-s})$, 
where $d_M$ denotes the distance to $M$ and $D^s$ is the total derivative of 
order $s\le r-1$ (Lemma 4.3 in [HW]). Then one solves $\dibar w_\e=u$ 
with $L^2$ estimate on $\omega_\e$ [H\"o1, H\"o2]. 
By the interior elliptic regularity of $\dibar$
this implies the uniform estimate $|w_\e|=o(\e^{r-m/2-1})=o(\e^k)$ 
on $\omega_{c\e}$ for any fixed $c\in (0,1)$ 
(the second display on p.\ 16 in [HW]). 
If $0 < s\le k$ we also have $|D^s w_\e|=o(\e^{k-s})$ on 
$\omega_{c'\e}$ for a fixed $c'\in (0,c)$
(see e.g.\ Lemma 3.2 and the proof of Proposition 2.3 in [FL]).
The function $f_\e=u-w_\e$ is holomorphic on $\omega_\e$;
as $\e\to 0$, $f_\e$ converges to $f$ uniformly on $K_0$ 
and in the $\cC^r$-sense on $M$. 
Since $\omega_\e$ is Runge in $V:=\omega_{\e_0}$ (see the
Remark preceding Theorem 3.2), we can approximate $f_\e$ by 
$F_\e\colon V\to \C$ and obtain the desired approximating sequence 
on a fixed open set $V\subset \C^n$. This proves the special case 
of Theorem 3.2.

For $Y=\C^N$ the result follows immediately by applying it componentwise. 
The case when $X$ is a Stein manifold (and $Y=\C^N$) reduces to the 
special case by embedding $X$ as a closed complex submanifold in some $\C^n$, 
or by applying the proof in [HW] with the Stein \nbd\ basis furnished 
by Theorem 3.1 above. To prove the general case of Theorem 3.2 we need the following.

\proclaim LEMMA 3.3:
Let $h\colon Z\to X$ be a \holo\ submersion of a complex manifold $Z$ 
to a complex manifold $X$. Assume that $S=K_0\cup M \subset X$ satisfies 
the hypotheses of Theorem 3.2. Let $U$ be an open set in $X$ containing 
$K_0$ and let $f\colon U\cup M\to Z$ be a section which 
is \holo\ on $U$ and smooth of class $\cC^1$ on $M$. Then  
$f(S)$ has a Stein \nbd\ basis in $Z$.

\demo Proof: We may assume that $U$ is Stein
and $K_0$ is $\cO(U)$-convex. The submanifold $f(M)\subset Z$  
is projected by $h$ bijectively onto the totally real submanifold 
$M\subset X$ and hence $f(M)$ is totally real in $Z$.
Since $f|_U$ is holomorphic, $f(U)$ is a closed complex submanifold 
of $Z|_U=h^{-1}(U)$ and hence by [De, S] it has an open Stein \nbd\ 
$\wt U\subset Z|_U$. For any compact $\cO(U)$-convex subset $K\ss U$ 
the set $f(K)$ is holomorphically convex in $f(U)$ and 
hence (since $f(U)$ is a closed complex submanifold of $\wt U$)
also in $\wt U$.  Applying this with $K=K_0$, 
and also with $K=(K_0\cup S)\cap N$ for some compact \nbd\ $N\subset U$ 
of $K_0$, we see that $f(S) \subset Z$ satisfies the hypothesis 
of Theorem 3.1 and hence it has a basis of Stein \nbd s. 
\smallskip

We continue with the proof of Theorem 3.2.
Assume that $f\colon U\cup M\to Z$ is a section of 
$h\colon Z\to X$ which is holomorphic in $U$ and of class 
$\cC^r$ on $M$ for some $r\ge 1$. 
Fix a Stein \nbd\ $\Omega\subset Z$ of $f(S)$ furnished 
by Lemma 3.3 and embed $\Omega$ as a closed complex submanifold
of a Euclidean space $\C^N$. There is an open \nbd\ $\wt\Omega \subset\C^N$
of $\Omega$ and a holomorphic retraction $\phi\colon \wt\Omega\to \Omega$
[DG]. We consider $f$ as a map into $\C^N$ via the embedding $\Omega\hra\C^N$. 
The special case of Theorem 3.2 gives a sequence 
of \holo\ maps $g_j \colon V\to \C^N$ $(j\in \N)$ in an open \nbd\  
$V\subset X$ of $S$ such that 
$\lim_{j\to\infty} g_j|_S = f|_S$ (the convergence is uniform on $K_0$
and in the $\cC^k$-sense on $M$). Let $V_j=\{x\in V\colon g_j(x)\in \wt\Omega\}$; 
this is a \nbd\ of $S$ for sufficiently large $j\in \N$.  
The sequence of maps $\wt f_j=\phi \circ g_j \colon V_j\to \Omega$ 
then satisfies the conclusion of Theorem 3.2 except that   
$\wt f_j$ need not be a section of $h$. This is corrected by 
projecting the point $\wt f_j(x)$ to the fiber $Z_x=h^{-1}(x)$ 
by the \holo\ retraction provided by the following

\proclaim LEMMA 3.4:
Let $\Omega$ be a Stein manifold and $h\colon \Omega\to V$ a
\holo\ submersion onto a complex manifold $V$. 
There are an open Stein set $W \subset V\times \Omega$ 
containing $\Gamma :=\{(x,z) \in V\times \Omega\colon  h(z)=x\}$
and a \holo\ retraction $\pi\colon W\to \Gamma$ such that
$\pi(x,z)=(x,\pi_2(x,z))$ for every $(x,z)\in W$.

It follows that $h(\pi_2(x,z))=x$, i.e., $\pi_2(x,\cdotp)$ 
is a holomorphic retraction of an open \nbd\ of the fiber 
$h^{-1}(x)$ in $\Omega$ onto $h^{-1}(x)$ for every fixed $x\in V$. 
Assuming Lemma 3.4 for a moment we set
$f_j(x)= \pi_2(x,\wt f_j(x))$ for $j=1,2,\ldots$; 
these are holomorphic sections of $h$ in small open 
\nbd s $V_j\subset X$ of $S$ satisfying Theorem 3.2.

The version of Theorem 3.2 with interpolation along a subvariety 
$X_0$ not intersecting $M=S\bs K_0$ can be satisfied by using an  
embedding $\Omega\hra \C^N$ and writing $f= f_0 + \sum_{l=1}^{l_0} h_l g_l$
in a \nbd\ $V\subset X$ of $S$, where $f_0$ and $g_l$ are maps 
$V\to\C^N$ such that $f_0$ is holomorphic in $V$,
$g_l \in \cC^r(V)$ is  \holo\ over a \nbd\ of $K_0$ (where 
$f$ is \holo), and the functions $h_1,\ldots,h_{l_0} \in \cO(V)$ 
vanish to order $s$ on $X_0$ and satisfy 
$X_0\cap V=\{x\in V\colon h_l(x)=0,\ 1\le l\le l_0\}$.
(See Lemma  8.1 in [FP2, p.\ 660] for the details.) 
Applying Theorem 3.2 with $Y=\C^N$
we approximate each $g_l$ on $S$ by a sequence of sections  
$g_{l,j}$ $(j=1,2,3,\ldots)$ holomorphic in a \nbd\ of $S$ in $X$. 
This gives a holomorphic sequence 
$\wt f_j = f_0  + \sum_{l=1}^{l_0} h_l g_{l,j}$ whose restriction
to $S$ converges to $f$ as $j\to \infty$ (in the sense of Theorem 3.2).
It remains to compose $\wt f_j$ with the two retractions as above
(first by the retraction $\phi$ onto $\Omega$ and then by the 
retraction furnished by Lemma 3.4) 
to get a sequence of \holo\ sections $f_j \colon V_j \to Z$
$(j=1,2,\ldots)$ which agree with $f$ to order $s$ 
on $X_0$ and satisfy $f_j|_S\to f|_S$ as $j\to \infty$. 
This completes the proof of Theorem 3.2 granted that Lemma 3.4 holds.

\demo Proof of Lemma 3.4: 
Clearly $\Gamma$ is a closed complex submanifold of $V\times\Omega$.
Let $p_1\colon V\times \Omega\to V$, 
$p_2 \colon V\times \Omega \to \Omega$ denote the 
respective projections. Observe that 
$p_2|_\Gamma \colon \Gamma\to \Omega$ is bijective
with the inverse $z\to (h(z),z)$; hence $\Gamma$ is 
Stein and consequently it has an open Stein \nbd\
in $V\times \Omega$. 

Consider the holomorphic vector subbundle 
$E \subset T(V\times\Omega)$ whose fiber over 
$(x,z)$ consists of all vectors $(0,\xi)$,
$\xi\in T_z \Omega$. The restricted bundle  
$E|_\Gamma$ can be decomposed as 
$E|_\Gamma= E_0\oplus N$ where $E_0=E|_\Gamma  \cap T\Gamma$
and $N$ is some complementary subbundle. Observe that
$N$ is the normal bundle of $\Gamma$ in $V\times \Omega$
(here we use the hypothesis that $h\colon \Omega\to V$
is a submersion). Let $N_0$ denote the zero section of $N$.
By the Docquier-Grauert theorem [DG] there is a biholomorphic 
map $\Phi$ from an open \nbd\ $U\subset N$ of $N_0$ onto an open 
\nbd\ $W\subset V\times \Omega$ of $\Gamma$ which maps 
$N_0$ onto $\Gamma$ and maps the fiber 
$N_{(x,z)} \cap U$ to $\{x\}\times \Omega =p_1^{-1}(x)$
for any $(x,z)\in \Gamma$. (We can obtain such $\Phi$ 
by composing local flows of holomorphic vector fields
which are tangent to $N$.) Choosing $U$ to have convex 
fibers it follows that $\Phi$ conjugates the base projection 
$U\to N_0$ to a \holo\ retraction $\pi\colon W\to \Gamma$ 
satisfying the conclusion of Lemma 3.4.
\smallskip

\it Remark. \rm
The loss of smoothness in Theorem 3.2 is due to the method 
of proof. When $K_0=\emptyset$ and $S$ is a \tr\ submanifold 
of $X$ of class $\cC^r$, the optimal $\cC^r$-to-$\cC^r$ approximation 
was proved by Range and Siu [RS] following the work of many authors. 
The kernel approach developed in [FL\O] gives precise approximation 
of smooth diffeomorphisms by biholomorphisms in tubular \nbd s. 
Both approaches can likely be adapted to the situation considered here, 
but we do not need such sharp approximation.

%
%
%
%
\beginsection 4. Proof of Theorem 2.1.

We first state parts (a) and (b) of Theorem 2.1 in a precise 
quantitative form; part (c) will be considered later.

\proclaim THEOREM 4.1: Let $X$ be a Stein manifold, $K\subset X$ a compact 
$\cO(X)$-convex subset and $Y$ a complex manifold 
with $\dim Y\le \dim X$. Choose a distance function $d$ on $Y$ induced by a 
complete Riemannian metric. Assume that $(f_0,\iota_0)\in \cS(X,Y)$ is such that
$f_0|_K \colon K\to Y$ is a holomorphic submersion and $\iota_0|_K=Tf_0|_K$. 
If $Y$ satisfies Property $S_n$ with $n=\dim X$ then for every $\e>0$ 
there is a homotopy $(f_t,\iota_t)\in \cS(X,Y)$ $(t\in [0,1])$ from 
$(f_0,\iota_0)$ to some $(f_1,Tf_1)\in \cS_{holo}(X,Y)$ such that 
for every $t\in [0,1]$, $f_t|_K \colon K\to Y$ is a \holo\ submersion, 
$Tf_t|_K=\iota_t|_K$, and $\sup_{x\in K} d(f_t(x),f_0(x)) <\e$.

\demo Proof: We shall follow the construction of \holo\ submersions $X\to\C^q$ 
in [F1], indicating the necessary changes and additional arguments.  

Assume that $U\subset X$ is an open set containing $K$ such that 
$f_0|_U\colon U\to Y$ is a \holo\ submersion and $Tf_0|_U=\iota_0|_U$. 
Choose a smooth \spsh\ Morse exhaustion function $\rho$ on $X$ 
such that $K\subset \{\rho<0\}\ss U$ and $0$ is a regular value of $\rho$.
We may assume furthermore that in some local holomorphic coordinates
near each critical point the function $\rho$ is a quadratic normal form
given by (6.1) in [F1].

The construction of the homotopy $(f_t,\iota_t) \in \cS(X,Y)$ 
is done in a countable sequence of {\it stages}, and each stage consists of 
finitely many {\it steps}. We use two different types of steps, one for crossing the 
noncritical levels of $\rho$ and the other one for crossing a critical level. 
The selection of sets involved in stages and steps is done in advance and 
depends only on the exhaustion function $\rho$.
On the other hand, some of the constants in the approximation at 
each stage (or step) are chosen inductively and depend on the partial 
solution obtained in the previous steps. The entire construction is quite 
similar to the proof of the Oka-Grauert principle in [HL2],
and especially in [FP1]. We first explain the global scheme;
compare with Sect.\ 6 in [F1].

Let $p_1,p_2,p_3,\ldots$ be the critical points of $\rho$ in $\{\rho>0\}\subset X$,
ordered so that $0<\rho(p_1)<\rho(p_2)<\rho(p_3)<\ldots$. Choose a sequence
$0 = c_0 < c_1 < c_2<\ldots$ with $\lim_{j\to\infty} c_j=+\infty$
such that $c_{2j-1} < \rho(p_j) < c_{2j}$ for every $j=1,2,\ldots$, and 
the numbers $c_{2j-1},\ c_{2j}$ are close to $\rho(p_j)$; the desired
closeness will be specified below when crossing the critical level 
$\{\rho=\rho(p_j)\}$. If there are only finitely many $p_j$'s, we choose 
the remainder of the sequence $c_j$ arbitrarily with $\lim_{j\to\infty} c_j=+\infty$. 
We subdivide the parameter interval 
of the homotopy into subintervals $I_j=[t_j,t_{j+1}]$ 
with $t_j = 1-2^{-j}$ $(j=0,1,2,\ldots)$.

In the $j$-th stage of the construction we assume inductively that we  
have a partial solution $(f_t,\iota_t)\in \cS(X,Y)$ for $t\in [0,t_j]$ 
satisfying 
$$
	\sup\{d(f_t(x),f_0(x)) \colon x\in K\} < \e_j, 
	\qquad t\in [0,t_j] 					\eqno(4.1)
$$ 
for some $\e_j \le \e (1-2^{-j-1})$, such that 
$f_{t_j} \colon \{\rho\le c_j\} \to Y$ is a \holo\ submersion 
and $\iota_{t_j}=Tf_{t_j}$ on $\{\rho\le c_j\}$. 
For $j=0$ these conditions are satisfied with $\e_0=\e/2$. 

Choose a number $\d_{j}\in (0,\e_j 2^{-j-1})$ such that 
any holomorphic map from $\{\rho\le c_j\}$ to $Y$ which is uniformly 
$\d_j$-close to $f_{t_j}$ (measured by the metric $d$ on $Y$) 
is a submersion on $\{\rho\le c_{j-1}\}$.  (For $j=0$ we choose $c_{-1}<0$
sufficiently  close to $0$ such that $K\subset \{\rho<c_{-1}\}$.) 
We also insure that $\d_j<\d_{j-1}/2$ where $\d_{j-1}$ is the analogous
number from the previous stage; this condition is vacuous for $j=0$. 
The goal of the $j$-th stage is to extend the solution 
$(f_t,\iota_t) \in \cS(X,Y)$ to the interval $t\in I_j=[t_j,t_{j+1}]$ 
(keeping it unchanged for $t\in [0,t_j]$) such that 
$$
	\sup\{ d(f_t(x),f_{t_j}(x)) \colon x\in X,\ \rho(x)\le c_j\} 
	< \d_{j}/2,  \qquad t\in I_j,  				    \eqno(4.2)
$$
$Tf_t=\iota_t$ on $\{\rho\le c_j\}$ for every $t\in I_j$, 
and $Tf_{t_{j+1}} =\iota_{t_{j+1}}$ on $\{\rho\le c_{j+1}\}$. 
Then $f_{t_{j+1}}$ satisfies the inductive hypothesis on 
the set $\{\rho\le c_{j+1}\}$ with  
$\e_{j+1} = \e_j+\d_j/2 < \e (1-2^{-j-2})$, 
thus completing the $j$-th stage. 

Assume for a moment that this process can be worked out.
For every $t \in [t_j,1)$ the 
map $f_t$ is \holo\ on $\{\rho\le c_j\}$ and $\iota_t=Tf_t$ there.
By (4.2) the limits $f_1=\lim_{t\to 1} f_t \colon X\to Y$ and 
$\iota_1=\lim_{t\to 1}\iota_t =\lim_{t\to 1} Tf_t \colon TX\to TY$  
exist uniformly on compacts in $X$. 
It follows that $f_1\colon X\to Y$ is holomorphic and $\iota_1=Tf_1$
on $X$. By the construction we also have 
$$ 
	d(f_1(x),f_{t_j}(x)) <\d_j,\qquad x\in \{\rho\le c_j\},\ j=0,1,2,\ldots.
$$
The choice of $\d_j$ implies that $f_1$ is 
a holomorphic submersion on $\{\rho\le c_{j-1}\}$. Since this holds 
for every $j\in \N$, $f_1$ is a holomorphic submersion of $X$ to $Y$
and hence $(f_1,Tf_1)\in \cS_{holo}(X,Y)$.
From (4.1) we also get $d(f_1(x),f_0(x))<\e$ for $x\in K$. 
This will complete the proof of Theorem 4.1
provided that we prove the inductive stage.

We first consider the {\it noncritical stages}, i.e., those 
for which $\rho$ has no critical values on $[c_j,c_{j+1}]$. 
(In our notation this happens for even values of $j$.) 
We solve the problem in finitely many steps of the following 
kind.  We have compact subsets $A,B$ of $X$ satisfying

\medskip
\item{(i)}  $A$, $B$, $C:=A\cap B$, and $\wt A:=A\cup B$ 
are (the closures of) strongly pseudoconvex domains in $X$,  
\item{(ii)}  $\bar{A\bs B} \cap \bar{B\bs A} =\emptyset$ (the separation property), and 
\item{(iii)} there is an open set $U$ in $X$ containing $B$ 
and a biholomorphic map $\psi\colon U\to U'$ onto an open 
subset $U'\subset \C^n$ containing the cube $Q$ (2.1) such that 
$\psi(B)$ is a convex subset of $Q$ and $O := \psi(A \cap U)\cap Q$ 
is a special convex set of the form (2.2).

\proclaim PROPOSITION 4.2: Let $A,B\subset X$ satisfy the  
properties (i)--(iii) above. Assume that $f\colon A\to Y$ is a holomorphic 
submersion. If $Y$ satisfies Property $S_n$ with $n=\dim X$ 
then for any $\d>0$ there is a homotopy of \holo\ submersions 
$f_t\colon A\to Y$ $(t\in [0,1])$, with $f_0=f$, such that 
$\sup_{x\in A} d(f_t(x),f(x))  <\d$ for every $t\in [0,1]$
and $f_1$ extends to a \holo\ submersion $A\cup B\to Y$.

For $Y=\C^q$ with $q<\dim X$ this is Lemma 6.3 in [F1].
Its proof also applies in our situation, except that we must replace 
the use of Proposition 3.3 in [F1] by the assumed Property $S_n$ 
of $Y$. We include a sketch of proof. 

By hypothesis the map $f\circ\psi^{-1} \colon U'\to Y$ is a holomorphic 
submersion on $O=\psi(A\cap U) \cap Q \subset \C^n$. By Property $S_n$ 
of $Y$ we can approximate it uniformly on some fixed \nbd\ of $O$ as close 
as desired by a \holo\ submersion $\wt g\colon Q\to Y$. (We actually apply $S_n$ 
with the pair $(rO,rQ)$ of dilated sets for some $r>1$ close to $1$ 
in order to get approximation on a \nbd\ of $O$.) Since $\psi(B) \subset Q$, 
the map $g=\wt g\circ\psi \colon B\to Y$ is a holomorphic submersion which 
approximates $f$ uniformly in a \nbd\ of $A\cap B$ as close as desired.  
By Lemma 5.1 in [F1] we have $f=g\circ \g$ for a biholomorphic map $\g$
defined in an open \nbd\ of $A\cap B$ in $X$ and uniformly very close
to identity map. (Its distance from the identity only depends on the 
distance of $f$ and $g$ on the given \nbd\ of $A\cap B$. 
The proof of Lemma 5.1 in [F1] holds for arbitrary
target manifold $Y$.) If the approximations are sufficiently close,
Theorem 4.1 in [F1] provides a decomposition $\g=\b\circ \a^{-1}$ 
in a \nbd\ of $A\cap B$, where $\a$ is a biholomorphic map close to the 
identity in a \nbd\ of $A$ in $X$ and $\b$ is a biholomorhic map close to the 
identity in a \nbd\ of $B$. Thus $f\circ\a=g\circ \b$ in a \nbd\ of $A\cap B$, 
and hence the two sides define a \holo\ submersion $\wt f \colon A\cup B \to\C^q$
which approximates $f$ on $A$. Furthermore, there is a homotopy 
$\a_t$ ($t\in [0,1]$) of biholomorphic maps close to the identity
in some fixed \nbd\ of $A$ such that $\a_0$ is the identity and
$\a_1=\a$. (It suffices to embed $X$ as a complex submanifold in some 
Euclidean space, take the convex linear combinations of $\a$ with the 
identity map, and project this homotopy back to the submanifold $X$ 
by a holomorphic retraction.) Then $f_t:= f\circ \a_t \colon A\to Y$ 
is a homotopy of holomorphic submersions from $f_0=f$ to $f_1=\wt f$ 
satisfying Proposition 4.2. This completes the proof.
\smallskip

It remains to explain how Proposition 4.2 is used in the $j$-th stage 
of the construction. Since $\rho$ is assumed to have no critical 
values in $[c_j,c_{j+1}]$, we can obtain the set $\{\rho\le c_{j+1}\}$ from 
$\{\rho\le c_j\}$ by finite number of attachings of `convex bumps'
of the type introduced just before Proposition 4.2. 
We begin with $A_0=\{\rho\le c_j\}$ and attach a bump $B_0$ 
to get $A_1=A_0\cup B_0$; then we attach a new bump $B_1$ to $A_1$ 
to obtain $A_2=A_1\cup B_1$, etc., until reaching 
$A_{k_j}=\{\rho\le c_{j+1}\}$ (see Lemma 12.3 in [HL1]). 
The required number of steps $k_j$ depends only on $\rho$.
We also subdivide the parameter interval 
$I_j=[t_j,t_{j+1}]$ into adjacent subintervals 
$I_{j,k}=[t_{j,k-1},t_{j,k}]$ $(k=1,2,\ldots,k_j)$
of equal length, one for every step. 

Assume inductively that for some $1\le k <k_j$ a solution 
$(f_t,\iota_t) \in \cS(X,Y)$ has already been defined for 
$t\in [0,t_{j,k-1}]$ such that $f:=f_{t_{j,k-1}}$ is a
holomorphic submersion from $A_{k-1}$ to $Y$ 
and $\iota:=\iota_{t_{j,k-1}} = Tf$ on this set.
Applying Proposition 4.2 on $A=A_{k-1}$ we extend the 
family of solutions to the next subinterval $t\in I_{j,k}$ such that 
$$
	\sup\{d( f_t(x), f_{t_{j,k-1}}(x) \colon x\in A_{k-1}\} < \d_j/2k_j
$$ 
(compare with (4.2)). We can also define $\iota_t$ for $t\in I_{j,k}$
such that $\iota_t=Tf_t$ on $A_{k-1}$, $\iota_{t_{j,k}}=Tf_{j,k}$ on $A_k$, 
$\iota_t$ is homotopic to $\iota_0$ and agrees with $\iota_0$ 
outside of some small \nbd\ of $A_k$. In $k_j$ steps we extend the 
family of solutions $(f_t,\iota_t)$ to $t\in I_j$ and 
thus complete the $j$-th stage.

It remains to consider the {\it critical stages}, i.e., those 
for which $\rho$ has a critical point $p$ with 
$c_j< \rho(p)< c_{j+1}$. For $Y=\C^q$ this is explained in 
sections 6.2--6.4 in [F1]. Since the proof needs a few modifications, we 
shall go through it step by step. 

Write $c=c_j$. We may assume that $c$ has been chosen as close to 
$\rho(p)$ as will be needed in the sequel.  Near $p$ we use local 
holomorphic coordinates on $X$ in which
$p=0$ and $\rho$ is a quadratic normal form (see (6.1) in [F1]). 
In particular, the stable manifold of $p=0$ for the gradient flow 
of $\rho$ is $\R^\nu \subset \C^n$ (the subspace spanned by the real 
parts of the first $\nu$ variables), where $\nu$ is the 
Morse index of $\rho$ at $p$. To cross the critical level 
$\rho=\rho(p)$ we perform the following three steps.

\medskip {\it Step 1: Extension to a handle.} 
We attach to $\{\rho\le c\}$ the disc $M\subset \R^\nu \subset \C^n$
(in the local coordinates) such that the attaching sphere $bM\subset \{\rho=c\}$ 
is complex tangential in the latter hypersurface (Subsect.\ 6.2
in [F1]). Let $(f,\iota)$ be the partial solution obtained after the
first $j-1$ stages, so $f$ is a \holo\ submersion from a \nbd\ of 
$\{\rho\le c\}$ to $Y$ and $\iota=Tf$ there.

\proclaim LEMMA 4.3: There is a \nbd\ $U$ of 
$S:=\{\rho\le c\} \cup M$ in $X$ and a smooth map $g \colon U\to Y$ 
which agrees with $f$ in a \nbd\ of $K_0:=\{\rho\le c\}$ such that 
for every $x\in M$ the differential $d g_x \colon T_x X\to T_{g(x)} Y$ 
is a surjective $\C$-linear map. Furthermore, $(f,\iota)$ can be connected
to $(g,Tg)$ by a path (homotopy) in $\cS(U,Y)$ which is fixed 
in a \nbd\ of $K_0$.

\demo Proof: The main point is the extension of $f$ and its 
$1$-jet to the handle $M$ such that the above properties are
satisfied. One uses Gromov's {\it convex integration lemma} ([Gr2]; 
Section 2.4 of [Gr3],  especially (D) and (E) in [Gr3, 2.4.1.];
Sect.\ 18.2 in [EM], especially Corollary 18.2.2.). 
The details given in [F1] (Lemmas 6.4 and 6.5)
for the case $Y=\C^q$ remain valid for arbitrary target 
manifold $Y$. The differential relation controlling the problem is 
{\it ample in the coordinate directions} 
and hence Gromov's lemma applies. 

It is (only) at this point
of the proof that we use the hypothesis on the existence of a 
fiberwise surjective map $\iota\colon TX\to TY$ with base
map $f$. When $q=\dim Y\le [{n+1\over 2}]$ it suffices to apply 
Thom's jet transversality theorem as in [F1], and in this case $\iota$
automatically exists.

\medskip {\it Step 2: Holomorphic approximation on a handlebody.}  
We denote the result of Step 1 again by $(f,\iota)$; 
thus $f$ is a \holo\ submersion from a \nbd\ of the strongly 
pseudoconvex domain $K_0=\{\rho \le c\} \subset X$ to $Y$, 
it is smooth in a \nbd\ of the handle $M$, $df_x\colon T_x X \to T_{f(x)} Y$ 
is surjective and $\C$-linear at each point $x\in M$, and $\iota_x=df_x$ 
for every $x\in S=K_0\cup M$.

By [F1, Subsect.\ 6.3] the pair $(K_0,M)$ 
satisfies the hypothesis of Theorem 3.2 above, and hence there is  
a holomorphic map $\wt f$ from a \nbd\ of $S$ to $Y$ which 
approximates $f$ uniformly on $K_0$ and in the $\cC^1$-sense on $M$ 
as close as desired. (In [F1] the corresponding approximation result 
was proved for $Y=\C^q$.) If the approximations are sufficiently close 
then  $\wt f$ is a holomorphic submersion from some open \nbd\ of $S$
to $Y$. We patch the new map with the old one outside of a \nbd\ of $S$.

\smallskip {\it Step 3: Crossing the critical level of $\rho$.} 
We denote the result of Step 2 again by $f$ and set $\iota=Tf$
in a \nbd\ of $S=\{\rho\le c\}\cup M$. We patch $\iota$ with $\iota_0$ 
outside of a \nbd\ of $S$ by using a cut-off function 
in the parameter of the homotopy connecting $\iota$ to $\iota_0$. 

By the scheme explained in [F1, Subsect.\ 6.4] we can approximately 
extend $f$ across the critical level of $\rho$ at $p$ by performing 
a noncritical stage with another \spsh\ function 
$\tau$ furnished by Lemma 6.7 in [F1]. Afterwards we revert back 
to the original exhaustion function $\rho$ and 
proceed to the next noncritical stage associated with $\rho$. 
When changing the domain of the solution 
(first from a \nbd\ of the handlebody $S$ to a suitable sublevel set of $\tau$, 
and later from a higher sublevel set of $\tau$ to a supercritical sublevel 
set $\rho$) we sacrifice a part of the domain, but the loss 
is compensated in the next stage. The approximation 
conditions for the given stage can be satisfied with the appropriate choices of 
constants at every step. No further changes from [F1] are needed
apart from those already explained for a noncritical stage. 
This completes the proof of Theorem 4.1 and hence of (a), (b) in 
Theorem 2.1. 

To prove part (c) of Theorem 2.1 we perform the same construction
for homotopies of submersions connecting a given pair 
$f_0,f_1\colon X\to Y$. Property $HS_n$ is needed to obtain 
the analogue of Proposition 5.2 for such homotopies. The remaining steps
in the proof, including the crossing of a critical level of $\rho$,
require only inessential modifications which we leave to the reader.

\medskip {\it Remarks.} 1. It is possible to prove a stronger version of 
Theorem 4.1 with interpolation of a given holomorphic submersion along 
a closed complex subvariety $X_0\subset X$; compare with Theorem 2.5 in [F1].  

\ni 2. Our proof also applies in the equidimensional case $\dim X=\dim Y=n$ 
as long as $Y$ satisfies Property $S_n$. Unfortunately we don't know any 
such example; the main case to be solved is $Y=\C^n$ 
(see Problems 1--3 in [F1]).

%
%
%
%
\beginsection 5. Approximation of submersions on subsets of $\C^n$.

In this section we prove Property $S_n$ and $HS_n$ 
for certain algebraic manifolds, in particular for those listed in 
Proposition 2.5. 

A complex algebraic subvariety $A$ in an algebraic manifold $V$ 
will be called {\it thin} if it does not contain any complex hypersurfaces, 
i.e., $A$ has complex codimension at least two in $V$. The following is a key lemma.

\proclaim LEMMA 5.1:  Assume that $O\subset Q\subset \C^n$ are as in 
(2.1), (2.2). If $A$ is a thin algebraic subvariety of $\C^n$
with $A\cap O=\emptyset$ then any holomorphic submersion 
$h\colon \C^n\bs A \to Y$ to a complex manifold $Y$ 
can be approximated uniformly on $O$ by holomorphic 
submersions $Q\to Y$ (and hence by submersions $\C^n\to Y$).

\demo Proof: Let $z=(z',z'')$ be the coordinates on $\C^n$ with  
$z'=(z_1,\ldots, z_{n-2})$ and $z''=(z_{n-1},z_n)$.
After a small linear change of coordinates the projection 
$\pi\colon \C^n\to \C^{n-2}$, $\pi(z',z'')=z'$, 
is proper when restricted to $A$. In this situation Lemma 3.4 from [F1] 
gives for any $\d>0$ a \holo\ automorphism $\psi$ of $\C^n$ of the form 
$\psi(z',z'')=(z',\b(z',z''))$ satisfying $\sup_{z\in O} |\psi(z)-z| <\d$
and $\psi(Q)\cap A =\emptyset$. Then $h\circ\psi \colon Q\to Y$ 
is a \holo\ submersion whose restriction to $O$ is uniformly 
close to $h$. This proves Lemma 5.1.
\medskip

Suppose now that $Y$ is a projective algebraic manifold of dimension 
$q$. Given a holomorphic submersion $f\colon O\to Y$ from a special 
convex set $O\subset \C^n$ for some $n>q$, our goal is to approximate
$f$ uniformly on $O$ by a rational map $h\colon \C^n \to Y$ 
which is a \holo\ submersion outside of a thin algebraic 
subvariety $A\subset \C^n$ not intersecting $O$. 
If such approximations exist then by Lemma 5.1 we can approximate
$h$ (and hence $f$) uniformly on $O$ by submersions $g\colon Q\to Y$,
thus proving that $Y$ satisfies Property $S_n$. In a similar way 
we establish $HS_n$.

The possibility of approximating $f$ by $h$ is of course 
a nontrivial condition on $Y$ which fails for example on 
Kobayashi-Eisenman hyperbolic manifolds.  
We shall first establish Property $S_n$ of the projective spaces 
$\CP^q$ when $n>q$; the proof for Grassmanians and their Zariski open subsets
with thin complements will follow the same pattern.  
(For $Y=\C^q$ see Proposition 3.3 in [F1].)

The quotient projection $\pi\colon \C^{q+1}_* = \C^{q+1}\bs \{0\} \to \CP^q$ 
is a holomorphic fiber bundle with fiber $\C^*$, and by adding the zero 
section we obtain the universal line bundle $L\to \CP^q$.
Assume that $f\colon O \to \CP^q$ is a holomorphic map 
on a compact convex set $O\subset \C^n$.
Since $O$ is contractible, the bundle $f^*L \to O$ is 
topologically trivial and hence holomorphically trivial [O, G3]. 
Therefore $f^*L$ admits a nowhere vanishing holomorphic section
which can be viewed as a holomorphic map 
$\wt f\colon O\to \C^{q+1}_*$ satisfying
$f=\pi\circ\wt f$ (a holomorphic lifting of $f$). 
We  approximate $\wt f$ uniformly on $O$ by a polynomial map 
$P\colon \C^n\to \C^{q+1}$ 
and take $h=\pi\circ P \colon \C^n\bs P^{-1}(0) \to \CP^q$;
by construction $h|_O$ approximates $f$.

It remains to show that for a generic choice of $P$ the map 
$h$ is a submersion outside of a thin subvariety in $\C^n$.
We say that $P$ is {\it transverse to $\pi$} at a point $z\in\C^n$
if $P(z)\ne 0$ and $dP_z$ is transverse to the fiber of $\pi$ through $P(z)$. 
Note that $h=\pi\circ P$ is a submersion in a \nbd\ of $z\in \C^n$ \iff\ 
$P$ is transverse to $\pi$ at $z$. Since $f\colon O\to\CP^q$
is a submersion and hence $\wt f$ is transverse to $\pi$ on $O$, 
we may assume (by choosing $P$ sufficiently uniformly close to $\wt f$ 
on a \nbd\ of $O$) that $P$ is also transverse to $\pi$ on $O$. 
To complete the proof it suffices to show that for a generic 
choice of $P$ the `singularity set'
$$
	\Sigma_P=P^{-1}(0) \cup \{z\in \C^n\bs P^{-1}(0)
	\colon dP_z \ {\rm is\ not\ transverse\ to\ }\pi\}
$$
is a thin algebraic subset of $\C^n$ provided $n>q$.

Before proceeding we consider also the case when $Y=G_{k,m}$
is the Grassman manifold consisting of $k$-dimensional complex 
subspaces of $\C^m$. We apply the above proof with the fibration 
$\pi\colon V_{k,m}\to Y$ where $V_{k,m}$ is 
the Stiefel variety of all $k$-frames in $\C^m$ acted upon by 
the group $GL_k(\C)$. We can identify $V_{k,m}$ with a Zariski open subset
in $\C^{km}$ with thin complement $B$ (which consists of $k\times m$ 
matrices with rank less than $k$), and $\pi$ defines an algebraic 
foliation $\cF$ of $\C^{km}$ which is nonsingular on $V_{k,m}$
and has homogeneous leaves biholomorphic to $GL_k(\C)$. 
By Grauert's main theorem from [G3] we can lift any holomorphic map 
$f\colon O\to Y$ from a compact convex set $O\subset\C^n$
to a holomorphic map $\wt f\colon O\to V_{k,m}$ such that 
$f=\pi\circ\wt f$ (the argument is essentially the same as above).
Next we approximate $\wt f$ by a polynomial 
map $P\colon \C^n\to \C^{km}$. We define the `singularity set' 
$\Sigma_P\subset \C^n$ of $P$ as above, 
except that we replace the origin in the target space by 
$B=\cF_{sing}$. To complete the proof it remains to show that,
for $n> \dim Y$, the set $\Sigma_P$ is thin in $\C^n$ for a generic 
choice of $P$. The same proof applies to Zariski open sets $\Omega$ 
with thin complements in $G_{k,m}$; it suffices to add to $B$
the (thin) $\pi$-preimage of the complement of $\Omega$.

To complete the proof of Proposition 2.5 (at least the part 
concerning the Property $S_n$) we need the following.
Fix a positive integer $N\in \N$ and denote by $V=\cP(n,m,N)$ 
the vector space of all polynomial maps $P\colon \C^n\to\C^m$ 
of degree at most $N$; note that $V$ may be identified with 
a Euclidean space.

\proclaim LEMMA 5.2: 
Let $\cF$ be an algebraic foliation of codimension $q$ on $\C^m$ 
with thin singular locus $\cF_{sing} \subset\C^m$.
Given a polynomial map $P\colon \C^n\to \C^m$ let $\Sigma_P$ 
consist of all points $z\in \C^n$ such that $P(z)\in \cF_{sing}$ or
$dP_z$ is not transverse to the leaf of $\cF$ through $P(z)$. 
If $n>q$ then the set $\Sigma_P$ is thin for all $P$
outside of a proper algebraic subvariety of $\cP(n,m,N)$.

The analogous result holds for any algebraic subsheaf $\cF$ 
(not necessarily integrable) of the tangent sheaf of $\C^m$ 
for which $\cF_{sing}$ (the locus of points where $\cF$ is not 
locally free) is thin. The proof is a standard application of 
transversality arguments [A, Fo] and is included only for completeness.

\demo Proof of Lemma 5.2:
Let $J\simeq \C^{n+m+nm}$ denote the manifold of one-jets of 
\holo\ maps $\C^n\to\C^m$. For each \holo\ map 
$f\colon U\subset\C^n\to \C^m$ the associated one-jet extension is  
$$ 
	z\in U \to j^1_z f = (z,f(z),\di_1 f(z),\ldots,\di_n f(z)) \in J 
$$
where $\di_j={\di \over \di \z_j}$. Denote by 
$p_1\colon J\to\C^n$ (resp.\ $p_2\colon J\to \C^m$) the source point 
(resp.\ the image point) projection.  

Given an algebraic foliation $\cF$ on $\C^m$ of codimension $q$ with 
singularity set $\cF_{sing}$, let $\Sigma_{\cF} \subset J$ denote 
the subset consisting of $p_2^{-1}(\cF_{sing})$ together with 
all one-jets over points $w\in \cF_{reg}=\C^m\bs \cF_{sing}$ which are 
not transverse to $T_w \cF$ (the tangent space to the leaf of $\cF$ at $w$). 
We claim that $\Sigma_{\cF}$ is a thin algebraic subset of $J$ 
provided that $\cF_{sing}$ is thin and $n>q$. 
Clearly $p^{-1}_2(\cF_{sing})$ is thin. 
Furthermore, each point $w\in \cF_{reg}$ is contained in a Zariski
open set $U\subset \C^m$ such that the restricted foliation $\cF|_U$ is 
defined by one-forms $\omega_1,\ldots,\omega_q$ on $\C^m$ with polynomial 
coefficients which are pointwise independent at $w$. A one-jet $j^1_z f$ 
with $f(z)=w$ fails to be transverse to $T_w\cF$ \iff\ the $q\times n$ matrix
with entries $\langle\omega_j(w), \di_k f(z)\rangle$ 
$(1\le j\le q,\ 1\le k\le n)$ has rank less than $q$.
A simple count shows that the subvariety consisting of all
such matrices has codimension $|n-q|+1$ which is $\ge 2$ when $n>q$.
This gives over each point $w\in \cF_{reg}$ at least two independent 
algebraic equations for $\Sigma_{\cF}$, thus proving our claim.

We identify $V=\cP(n,m,N)$ with the Euclidean space whose elements 
are collections $c= \{ (c_\a)\colon |\a|\le N\}$, where $c_\a\in \C^m$ 
for each multiindex $\a=(\a_1,\ldots,\a_n)$; 
the correspondence is given by
$$ 
	c=\{c_\a \colon |\a|\le N\} \longrightarrow 
	P_c(z)=\sum_{|\a|\le N} c_\a\, z^\a \in V.
$$
Let $\gamma \colon \C^n\times V \to J$ associate to 
every pair $(z,P)$ the one-jet $j^1_z P\in J$.
Clearly $\gamma$ is polynomial in $z$ and linear in the coefficients of $P$. 
Furthermore, fixing $z\in \C^n$ and the coefficients $c_{\a}$ with 
$2\le |\a|\le N$, $\gamma$ gives a linear map of maximal rank 
(a linear submersion) from the space of coefficients $\{c_\a \colon |\a|\le 1\}$ 
of order $\le 1$ onto the fiber $p_1^{-1}(z)\subset J$. 
It follows that $\Sigma^1 \colon = \gamma^{-1}(\Sigma_{\cF})$ 
is a thin algebraic subvariety of $\C^n\times V$. 
Hence for every $P$ outside a proper algebraic subvariety of $V$
the set $\Sigma_P= \{z\in \C^n\colon (z,P)\in \Sigma^1\}$
is a thin subvariety of $\C^n$. By construction, $\Sigma_P$ is the set 
of point $z\in \C^n$ at which $dP_z$ fails to be transverse to the 
foliation $\cF$ (or $P(z)$ belongs to $\cF_{sing}$). This proves Lemma 5.2.

\smallskip
This establishes Property $S_n$ of any complex manifold $Y$ in 
Proposition 2.5 for all $n>\dim Y$. 
The proof that $\CP^1$ also satisfies Property $S_1$ requires a 
slightly different argument as follows. Let $f\colon O\to \CP^1=\C\cup\{\infty\}$ 
be a holomorphic submersion from a convex set $O\subset \C$ to the Riemann sphere.
Then $f^{-1}(\infty) \subset O$ is a finite subset of $O$ which we
may assume to be contained in the interior of $O$. Let $B\subset \C$ 
be a convex set such that $B\cap f^{-1}(\infty)=\emptyset$,
$B\cup O$ is also convex, and $\bar{B\bs O}\cap \bar{O\bs B}=\emptyset$. 
Let $C=B\cap O$. By Theorem 3.1 in [F1] we can approximate 
$f|_C \colon C\to \C=\CP^1\bs \{\infty\}$ by entire noncritical 
functions $g\colon \C\to\C$. Proceeding as in the proof of 
Proposition 4.2 above we find a biholomorphic transition map $\gamma$ 
between $f$ and $g$ on a \nbd\ of $C$, decompose $\gamma=\beta\circ\alpha^{-1}$
and thus patch $f$ and $g$ into a submersion $\wt f\colon B\cup O\to\CP^1$
which approximates $f$ on $O$. Furthermore we can arrange that 
$\wt f^{-1}(\infty)=f^{-1}(\infty)$.  In finitely many steps of this kind
we approximate $f$ by a submersion $Q\to \CP^1$ as desired.

We now show that every manifold $Y$ as above also satisfies Property $HS_n$
for all $n>q=\dim Y$. Recall that we have a submersion
$\pi\colon \C^m\bs B\to Y$ onto $Y$ where $B$ is thin in $\C^m$
and contains $\cF_{sing}$.
Let $O\subset Q\subset \C^n$ be as in (2.1), (2.2). Assume that 
$f_t\colon O\subset \C^n\to Y$ $(t\in [0,1])$ is  a homotopy of holomorphic 
submersions such that $f_0$ and $f_1$ extend to submersions $Q\to Y$. 
By the same argument as in the proof of Property $S_n$
we can lift $\{f_t\}$ to a homotopy of maps $\wt f_t\colon O\to \C^m$ 
$(t\in [0,1])$ which are transverse to the foliation $\cF$ defined by $\pi$
and such that $\wt f_0$, $\wt f_1$ are defined and transverse to $\cF$ on $Q$. 
We approximate  $\{\wt f_t\}$ by a homotopy of polynomial maps $P_t\colon \C^n\to \C^m$ 
which also depends polynomially on a parameter $t\in \C$ (the approximation 
of $\wt f_t$ by $P_t$ takes place on $O$ for $t\in (0,1)$ and on $Q$ 
for $t=0$ and $t=1$).   

Let $N$ be the maximal degree of the maps $P_t$ in the above
family, considered as polynomials both in $t\in \C$ and $z\in\C^n$.  
Denote by $\wt \cP(n,m,N)$ the space of all polynomial maps
$(t,z)\in \C\times \C^n \to\C^m$ of degree $\le N$. 
For every $P\in \wt \cP(n,m,N)$ and $t\in\C$ we have  
$P_t=P(t,\cdotp)\in \cP(n,m,N)$. The proof of Lemma 5.2 shows that 
for a generic choice of $P \in \wt \cP(n,m,N)$ the 
`singularity set' $\Sigma_P \subset \C^{1+n}$, consisting of all 
$(t,z)$ for which the one-jet $j^1_z P_t$ is not transverse 
to $\cF$ (or $P_t(z) \in B$), is a thin algebraic 
subvariety of $\C^{1+n}$. Hence for all but finitely many
$t\in \C$ the singularity set $\Sigma_t \subset \C^n$ of 
the map $P_t$ is also thin.
By a small deformation of the segment $[0,1] \subset \R$
inside $\C$ we may can avoid this finite exceptional set of $t$'s,
thus obtaining a homotopy $P_t\colon \C^n\to\C^m$ of polynomial
maps approximating $\wt f_t$ $(t\in [0,1])$ 
such that $\Sigma_t$ is thin for all $t\in [0,1]$. 
Hence $h_t=\pi\circ P_t \colon \C^n\bs \Sigma_t\to Y$
is a homotopy of submersions which approximates the original
homotopy $f_t \colon O\to Y$ uniformly on $O$, and uniformly
with respect to the parameter $t$. 

We can now conclude the proof as in Lemma 5.1: Applying [F1, Lemma 3.4] 
with the additional parameter $t\in [0,1]$ we obtain a family 
of holomorphic automorphisms $\psi_t(z',z'')=(z',\beta_t(z',z''))$
of $\C^n$, depending smoothly on $t\in [0,1]$, such that for every 
$t$ we have $\psi_t(Q)\cap \Sigma_t =\emptyset$,
$\psi_t|_O$ is close to the identity map on $O$, and the maps 
$\psi_0$ and $\psi_1$ are close to the identity map on $Q$. 
The homotopy of holomorphic submersions $h_t \circ \psi_t \colon Q\to  Y$ 
satisfy all the required properties, except perhaps the interpolation 
condition at the endpoints $t=0,1$ which is easily fixed.

When proving the Property $S_n$ or $HS_n$ for $\C^q$ 
(or a Zariski open subset with thin complement in $\C^q$) we proceed
as above but skip the first step, i.e, we can directly deal with
polynomial maps $\C^n\to\C^q$. This completes the proof of Proposition 2.5.

\smallskip
\it Remark. \rm 
Our proof of Property $HS_n$ breaks down for multi-para\-me\-ter
families of submersions: we still prove that for a generic choice of the polynomial
map $P\colon \C^k\times \C^n\to \C^m$ its singularity set 
$\Sigma_P \subset \C^k_t\times \C^n_z$ (with $t\in \C^k$ being the
parameter) is thin; hence for all $t$ outside of a proper algebraic subvariety 
$Z\subset \C^k$ the singularity set $\Sigma_t \subset\C^n$
of $P_t=P(t,\cdotp) \colon \C^n\to\C^m$ is thin in $\C^n$. However, 
when $k>1$ we may not be able to avoid the exceptional set $Z$ by a 
small deformation of the parameter cube $[0,1]^k\subset \R^k$ in $\C^k$.
Such a multi-parameter analogue of Property $HS_n$ would be needed to 
obtain the complete parametric homotopy principle for submersions $X\to Y$.

\beginsection 6. A simple proof of the Oka-Grauert-Gromov theorem.

We give a simple proof of the following result from [F3].

\proclaim THEOREM 6.1:
Let $h\colon Z\to X$ be a holomorphic fiber bundle over a Stein
manifold $X$. If the fiber $Y=h^{-1}(x)$ $(x\in X)$ admits a finite 
dominating collection of sprays then the inclusion 
${\rm Holo}(X,Z) \hra {\rm Cont}(X,Z)$ of the space
of holomorphic sections into the space of continuous
sections is a weak homotopy equivalence.

The spaces of sections are equipped with the compact-open topology.
For the definition of `dominating families of sprays' see [F3]. 
The classical case when $Z$ is a principal holomorphic bundle 
(with fiber a complex Lie group or homogeneous space) is due to 
Grauert [G3]. The case when the fiber $Y$ admits a 
dominating spray is due to Gromov [Gr4, Sec.\ 2.8], 
and a detailed proof can be found in [FP1].

\demo Proof: 
We use the scheme of proof of Theorem 2.1 in Sect.\ 4
(which is similar to the one in [FP1]).
We may assume that the fiber $Y$ is connected.

At a {\it noncriticial step} we have a continuous section $f\colon X\to Z$
which is holomorphic in an open \nbd\ of a smoothly 
bounded, compact, \spsc\ domain $A\subset X$. We attach 
a small `convex bump' $B$ satisfying
\item{(i)} $\bar{A\bs B}\cap \bar{B\bs A}=\emptyset$ and the union 
$A\cup B$ is again smoothly bounded \spsc,  
\item{(ii)} in some local holomorphic coordinates in a \nbd\ $U\subset X$ 
of $B$ the sets $B$ and $C=A\cap B$ are compact convex subset of $\C^n$ 
($n=\dim X$), and 
\item{(iii)} the restricted bundle $Z|_U$ is trivial, $Z|_U \simeq U\times Y$.

The set $B$ can be a `special convex set' of the form (2.2) in some local
coordinates on $X$. The goal is to approximate $f$ uniformly on $A$ 
as close as desired by a continuous section $\wt f \colon X\to Z$ which is 
holomorphic in a \nbd\ of $A\cup B$. This is accomplished 
by the `noncritical case' in [Gr4] or [FP1, Sect.\ 6] and here
we do not propose any changes. We recall the main steps for the 
sake of the reader. 

Let $f_1\colon U\to Z|_U\simeq U\times Y$ be any \holo\ 
section (for instance, a constant section). Since $C$ is 
convex and hence holomorphically contractible, there is a homotopy 
$\{f_t\}_{t \in [0,1]}$ of holomorphic sections over a \nbd\ $V\subset U$ of $C$ 
connecting $f_0:= f|_V$ to $f_1|_V$. Since $Y$ admits a finite dominating 
family of sprays, the homotopy version of the Oka-Weil theorem 
[F3, Theorem 3.1] gives a uniform approximation of the  homotopy $\{f_t\}$ 
on a smaller \nbd\ of $C$ by a homotopy of sections $\{\wt f_t\}_{t\in [0,1]}$ 
which are holomorphic over a \nbd\ of $B$, with $\wt f_1=f_1$
and $g:=\wt f_0$ very close to $f$ on a \nbd\ of $C$. 
(This approximation result is essentially due to Grauert [G1, G2]; 
see also [Gr4] and [FP1].)
Applying Theorem 4.1 in [F3] (or Theorem 5.1 in [FP1] when $Y$ admits 
a dominating spray) we glue the pair of sections $f$, $g$ 
into a section $\wt f$ which is \holo\ in an open \nbd\ of $A\cup B$
and extends to a continuous section over $X$. This completes the noncritical step.

A complication in this process arises when crossing a critical level
$\{\rho=c_0\}$ of the given \spsh\ exhaustion function $\rho\colon X\to \R$. 
In that case we have $A=\{\rho\le c\}$ for some $c< c_0$ close 
to $c_0$ (such that $\rho$ has no critical values on $[c,c_0)$).
The set $B$ is a `handle' attached to $A$ such that the attaching 
set $C=A\cap B$ is no longer contractible. (The union $A\cup B$ is 
diffeomorphic to the sublevel set $\{\rho\le c'\}$ for some $c'>c_0$.)
In this case we cannot find the desired homotopy $\{f_t\}$ as above. 
In [FP1] the difficulty was avoided by applying the noncritical case 
with an additional parameter to construct $\{f_t\}$, beginning at the 
`core' of $C$ (which is a \tr\ torus) and performing `approximation and gluing' 
until reaching $C$ in finitely many steps (Theorem 4.5 in [FP1]). 
An alternative method, proposed by Gromov [Gr4] and developed in [FP2], 
uses a more complicated induction scheme and remains applicable 
even if the submersion $Z\to X$ is not locally trivial, as long 
it admits fiber-dominating sprays over small open subsets of $X$.

Here we propose a simple alternative way to pass the critical level of $\rho$
by applying Theorem 3.2. Assume that $A$ and $f$ are as above.
We attach to $A$ a smooth \tr\ handle $M$ passing through 
the critical point $p$,
with $\dim M$ equal the Morse index of $p$ (see Sect.\ 4, 
{\it Step 1: extension to a handle}). By Theorem 3.2 we can 
approximate $f$ uniformly on $A\cup M$ by a section 
$\wt f$ which is \holo\ in a \nbd\ of $A\cup M$ and continuous on $X$ 
(compare with {\it Step 2} in Sect.\ 4). The proof is completed 
by {\it Step 3} in Sect.\ 4 without any changes 
(i.e., we use the noncritical case with a different \spsh\ function 
$\tau$ in order to pass the critical level of $\rho$ at $p$,  
then we revert back to $\rho$ and proceed by the noncritical case 
to the next critical level of $\rho$.) All steps adapt  easily to the 
parametric case and hence we obtain the full statement of Theorem 1.4 in [FP1]
under the weaker assumption that the fiber admits a dominating family of 
sprays (instead of a dominating spray).  This completes the proof.
\smallskip

Comparing with the proof of Theorem 2.1
we see that the only essential difference lies in the method 
of local approximation and gluing of pairs of sections. 
In the proof of Theorem 6.1 we use a dominating family of 
sprays on the fiber to linearize the problem. 
On the other hand, in Theorem 2.1 we patch 
a pair of submersions by decomposing the transition map between them
in the source manifold $X$ without using any properties
of the target manifold (these are only used for the approximation).  
The method of passing a critical level is identical in both proofs, 
except that the construction of submersions requires 
a maximal rank extension of the one-jet across the handle. 

The above may be conceptually the simplest available 
proof of the Oka-Grauert-Gromov principle, even in Grauert's classical 
case concerning principal bundles with homogeneous fibers. 
Unlike in the earlier papers, holomorphic sections
are constructed without using the intermediate technical devices 
for multi-parameter homotopies of sections.

\smallskip
{\it Acknowledgements.}
This research has been supported in part by a grant from the 
Ministry of Science and Education of the Republic of Slovenia. 
Kind thanks to J.\ Winkelmann for helpful discussions regarding 
Property $S_n$.

%
%
%
%
\medskip
\def\same{{--------\ }}

\beginsection References

\item{[A]} ABRAHAM, R., 
Transversality in manifolds of mappings. 
{\it Bull.\ Amer.\ Math.\ Soc.}, 69 (1963), 470--474.

\item{[B]} BARTH, W., PETERS, C., VAN DE VEN, A.,
{\it Compact Complex Surfaces.}
Springer, Berlin--Heidelberg--New Zork--Tokyo, 1984. 

\item{[De]} DEMAILLY, J.-P.,
Cohomology of $q$-convex spaces in top degrees.
{\it Math.\ Z.}, 204 (1990), 283--295.

\item{[DG]} 
DOCQUIER, F., GRAUERT, H.,
Levisches Problem und Rungescher Satz f\"ur 
Teilgebiete Steinscher Mannigfaltigkeiten. (German) 
{\it Math.\ Ann.}, 140 (1960), 94--123.

\item{[E]} ELIASHBERG, Y.,
Topological characterization of Stein manifolds of dimension $>2$.
{\it Internat.\ J.\ Math.}, 1 (1990), 29--46.

\item{[EM]}
ELIASHBERG, Y., MISHACHEV, N., 
{\it Introduction to the $h$-principle.} 
Graduate Studies in Math., 48.
Amer.\ Math.\ Soc., Providence, RI, 2002.

\item{[Fo]} 
FORSTER, O.,
Plongements des vari\'et\'es de Stein.  
{\it Comment.\ Math.\ Helv.}, 45 (1970), 170--184.

\item{[F1]} 
FORSTNERI\v C, F., 
Noncritical holomorphic functions on Stein manifolds.
{\it Acta Math.}, 191 (2003). 
[arXiv: math.CV/0211112]

\item{[F2]}  \same,
The homotopy principle in complex analysis: A survey, in 
\it Explorations in Complex and Riemannian Geometry: A Volume 
dedicated to Robert E.\ Greene \rm 
(J.\ Bland, K.-T.\ Kim, and S.\ G.\ Krantz, eds., 73--99).
Contemporary Mathematics, 332, 
American Mathematical Society, Providence, 2003.

\item{[F3]} \same,
The Oka principle for sections of subelliptic submersions.
{\it Math.\ Z.}, 241 (2002), 527--551.

\item{[FK]} FORSTNERI\v C, F., KOZAK, J.,
Strongly pseudoconvex handlebodies.
{\it J.\ Korean Math.\ Soc.}, 40 (2003), 727--746.

\item{[FL\O]}
FORSTNERI\v C, F., L\O W, E., \O VRELID, N.,
Solving the $d$- and $\overline\partial$-equa\-tions in thin tubes 
and applications to mappings. 
{\it Michigan Math.\ J.}, 49 (2001), 369--416. 

\item{[FP1]} 
FORSTNERI\v C, F., PREZELJ, J., 
Oka's principle for holomorphic fiber bundles with sprays.
{\it Math.\ Ann.}, 317 (2000), 117-154.

\item{[FP2]} \same,  
Oka's principle for holomorphic submersions with sprays.
{\it Math.\ Ann.}, 322  (2002), 633-666.

\item{[FP3]} 
\same, 
Extending holomorphic sections from complex subvarieties.
{\it Math.\ Z.}, 236  (2001), 43--68.

\item{[FR]} 
FORSTNERI\v C, F., ROSAY, J.-P.,
Approximation of biholomorphic mappings by automorphisms of $\C^n$.
{\it Invent.\ Math.}, 112  (1993), 323--349.
Erratum, {\it Invent.\ Math.}, 118  (1994), 573--574.

\item{[G1]} GRAUERT, H., 
Approximationss\"atze f\"ur holomorphe Funktionen mit \break
Werten in komplexen R\"aumen.
{\it Math.\ Ann.}, 133 (1957), 139--159.

\item{[G2]} 
\same, 
Holomorphe Funktionen mit Werten in komplexen Lieschen Gruppen.
{\it Math.\ Ann.}, 133 (1957), 450--472.

\item{[G3]} 
\same,
Analytische Faserungen \"uber holomorph-vollst\"andigen R\"aumen.
{\it Math.\ Ann.}, 135 (1958), 263--273.

\item{[Gr1]} 
GROMOV, M.,
Stable maps of foliations into manifolds.
{\it Izv.\ Akad.\ Nauk, S.S.S.R.}, 33 (1969), 707--734.

\item{[Gr2]} 
\same,
Convex integration of differential relations, I. (Russian) 
{\it Izv.\ Akad.\ Nauk SSSR Ser. Mat.}, 37 (1973), 329--343. 
English translation in {\it Math.\ USSR--Izv.}, 7 (1973), 329--343.

\item{[Gr3]} 
\same,  
{\it Partial Differential Relations.}
Ergebnisse der Mathematik und ihrer Grenzgebiete (3), 9.
Springer, Berlin--New York, 1986.

\item{[Gr4]} 
\same, 
Oka's principle for holomorphic sections of elliptic bundles.
{\it J.\ Amer.\ Math.\ Soc.}, 2 (1989), 851-897.

\item{[GN]} 
GUNNING, R.\ C., NARASIMHAN, R., 
Immersion of open Riemann surfaces. 
{\it Math.\ Ann.}, 174  (1967), 103--108. 

\item{[GR]} 
GUNNING, R.\ C., ROSSI, H., 
{\it Analytic functions of several complex variables.}
Prentice--Hall, Englewood Cliffs, 1965.

\item{[HL1]} 
HENKIN, G.\ M., LEITERER, J.,
{\it Andreotti-Grauert Theory by Integral Formulas.}
Progress in Math., 74, Birkh\"auser, Boston, 1988.

\item{[HL2]} \same, 
The Oka-Grauert principle without induction over the basis dimension.
{\it Math.\ Ann.}, 311 (1998), 71--93.

\item{[H\"o1]} 
H\"ORMANDER, L., 
$L\sp{2}$ estimates and existence theorems for the $\bar \partial$ operator. 
{\it Acta Math.}, 113 (1965), 89--152. 

\item{[H\"o2]} \same, 
{\it An Introduction to Complex Analysis in Several Variables}. 
Third ed. 
North Holland, Amsterdam, 1990.

\item{[HW]} 
H\"ORMANDER, L., WERMER, J.,
Uniform approximations on compact sets in $\C^n$.
{\it Math.\ Scand.}, 23 (1968), 5--21.

\item{[O]} OKA, K., 
Sur les fonctions des plusieurs variables. 
III: Deuxi\`eme probl\`eme de Cousin.
{\it J.\ Sc.\ Hiroshima Univ.}, 9 (1939), 7--19.

\item{[P]} 
PHILLIPS, A.,
Submersions of open manifolds.
{\it Topology}, 6 (1967), 170--206.
 
\item{[RS]}  
RANGE, R.\ M., SIU, Y.\ T., 
$\cC^k$ approximation by holomorphic functions and $\bar \partial$-closed 
forms on $\cC^k$ submanifolds of a complex manifold.
{\it Math.\ Ann.}, 210 (1974), 105--122.

\item{[Ro]} 
ROSAY, J.-P.,
A counterexample related to Hartog's phenomenon
(a question by E.\ Chirka). 
{\it Michigan Math.\ J.}, 45 (1998), 529--535.

\item{[S]} SIU, J.-T.,
Every Stein subvariety admits a Stein neighborhood.
{\it Invent.\ Math.}, 38 (1976), 89--100.

\item{[W]} WINKELMANN, J.,
The Oka-principle for mappings between Riemann surfaces.
{\it Enseign.\ Math.} (2) 39 (1993), 143--151.

\bigskip\medskip
{\it Address:}
Institute of Mathematics, Physics and Mechanics,
University of Ljubljana, Jadranska 19, 1000 Ljubljana, Slovenia

\medskip
{\it E-mail:} franc.forstneric@fmf.uni-lj.si

\bye